\documentclass[conference]{IEEEtran}

\usepackage[utf8]{inputenc}
\usepackage{amsmath}
\usepackage{amsfonts}
\usepackage{amssymb}
\usepackage{textcomp}
\usepackage{xcolor}
\usepackage{caption}
\usepackage{subcaption}
\usepackage{gensymb}
\usepackage{algorithm2e}
\usepackage{graphicx}
\usepackage{lipsum}
\graphicspath{ {./images/} }

\ifCLASSOPTIONcompsoc

\else

\fi
\usepackage[backend=biber,sorting=none,style=ieee,url=false]{biblatex}
\usepackage[utf8]{inputenc}
\usepackage[english]{babel}

\usepackage[utf8]{inputenc}
 \usepackage{fancyhdr}
\begin{document}

\makeatletter
\newcommand{\linebreakand}{%
  \end{@IEEEauthorhalign}
  \hfill\mbox{}\par
  \mbox{}\hfill\begin{@IEEEauthorhalign}
}
\makeatother

\author{
\IEEEauthorblockN{Tanvi Agrawal\IEEEauthorrefmark{1}}
    \IEEEauthorblockA{\textit{Research Head - thrustMIT}\\
    \textit{Department of Mechanical} \\
    \textit{and Manufacturing Engineering}\\
    \and
    \IEEEauthorblockN{Utkarsh Anand\IEEEauthorrefmark{2}}
    \IEEEauthorblockA{\textit{Payload Head - thrustMIT}\\
    \textit{Department of Electrical} \\
    \textit{and Electronics Engineering}\\}
    \linebreakand
    \IEEEauthorblockN{Dr. Srinivas G\IEEEauthorrefmark{0}}
    \IEEEauthorblockA{\textit{Faculty Advisor - thrustMIT}\\
    \textit{Asst. Prof- Sr. Scale} \\
    \textit{Department of Aeronatical}\\
    \textit{and Automobile Engineering}\\}
    }

    \linebreakand
    \\
    Manipal Institute of Technology,\\Manipal Academy of Higher Education, Manipal, Karnataka - 576104, India\\

    Email: 
    \IEEEauthorrefmark{1}tanviagrawal02@gmail.com
    \IEEEauthorrefmark{2}utkarshanand221@gmail.com,
    \linebreakand
    Corresponding Author email: srinivas.g@manipal.edu

}
\title{Optimization of a Runge-Kutta 4th Order Method-based Airbrake Control System for High-Speed Vehicles Using Neural Networks}
\maketitle

\begin{abstract}
The Runge-Kutta 4th Order (RK4) technique is extensively employed in the numerical solution of differential equations for airbrake control system design. However, its computational efficacy may encounter restrictions when dealing with high-speed vehicles that experience intricate aerodynamic forces. Using a Neural Network, a unique technique to improving the RK4-based airbrakes code is provided. The Neural Network is trained on numerous aspects of the high-speed vehicle as well as the current status of the airbrakes. This data was generated through the traditional RK4-based simulations and can predict the state of the airbrakes for any given state of the rocket in real-time. The proposed approach is demonstrated on a high-speed airbrakes control system, achieving comparable or better performance than the traditional RK4-based system while significantly reducing computational time by reducing the number of mathematical operations. The proposed method can adapt to changes in flow conditions and optimize the airbrakes system in real-time.
{ Keywords: Airbrakes, Control System, High-speed vehicle, Runge-Kutta 4th Order, Neural Network}

\end{abstract}

\section{Introduction}

Airbrakes are essential safety mechanisms for high-speed vehicles, such as aircraft and sounding rockets. They help to control and reduce the speed of the vehicle by converting its kinetic energy into other forms of energy. Airbrakes control systems rely on solving complex differential equations that describe the physics of the vehicle's motion and aerodynamic forces. The Runge-Kutta 4th Order (RK4) method is a prevalent numerical approach utilized in the resolution of these differential equations within the context of airbrake control system design\cite{arif1996optimization}. However, the computational efficiency of the RK4 method can be limited, especially for high-speed vehicles with complex aerodynamic forces.\cite{haykin1992solving}

In recent times, there has been a remarkable surge in the field of Artificial Intelligence, with Neural Networks emerging as formidable instruments for addressing intricate problems\cite{lecun1989optimal}.The objective of this study is to furnish empirical substantiation concerning the effectiveness of the proposed methodology in governing high-speed airbrakes. We will conduct a comparative analysis of the performance of a neural network-based airbrakes control system with that of a conventional RK4-based system, and evaluate the findings. Moreover, we will investigate the adaptability of the proposed technique to varying flow conditions, as well as its ability to optimize the airbrakes system in real-time.

This research proposes a novel and effective method for optimizing airbrakes control systems that utilize the RK4 method by incorporating neural networks. This method has the potential to reduce computational time and improve the accuracy of airbrakes control systems for high-speed vehicles. Furthermore, it can be implemented for other systems requiring real-time control and optimization.

\section{Background Theory}

\subsection{Runge-Kutta 4th Order Method}

The Runge-Kutta 4th Order (RK4) method is a numerical approach employed in various scientific and engineering domains, such as physics, chemistry, and control systems. It is particularly valuable for solving ordinary differential equations (ODEs) when analytical solutions are either unavailable or challenging to compute. The RK4 method involves generating four intermediate estimates of the dependent variable, which are subsequently weighted and combined to obtain the final estimate at the subsequent time step \cite{BUTCHER1996113}.

Consider a dependent variable, y(t), representing the value at time t. The RK4 method determines the value of y at the next time step, t + h, through a series of equations. Initially, k1 is computed by evaluating the product of the step size, h, and the rate of change function, f(t, y), at the current time and state \cite{doi:10.1080/00207169108803994}. 
The RK4 method computes the value of y at the next time step, t + h, as follows:
\begin{equation}
   k1 = h * f(t, y) 
\end{equation}
\begin{equation}
    k2 = h * f(t + h/2, y + k1/2)
\end{equation}
\begin{equation}
    k3 = h * f(t + h/2, y + k2/2)
\end{equation}
\begin{equation}
   k4 = h * f(t + h, y + k3) 
\end{equation}

\begin{equation}
   y(t + h) = y(t) + (k1 + 2k2 + 2k3 + k4) / 6 
\end{equation}

where f(t, y) is the function that describes the rate of change of y with respect to t, and h is the step size. The intermediate estimates k1, k2, k3, and k4 are calculated at different points within the time step, and their weighted sum is used to update the value of y at the next time step\cite{MUNTHEKAAS1999115}

Notably, the RK4 method exhibits fourth-order precision, implying that the error in the approximation is directly proportional to the fourth power of the step size, h. This characteristic renders the RK4 method highly accurate and efficient, particularly for systems characterized by complex dynamics\cite{doi:10.1080/00207169508804437}.

The widespread use of the RK4 method in scientific and engineering disciplines underscores its significance. Researchers extensively rely on the RK4 method to simulate and analyze intricate systems that lack analytical solutions or present challenges in their calculation\cite{doi:10.1080/00207169908804817}.

Within the domain of airbrake control systems for high-speed vehicles, the RK4 method serves as a valuable tool for modeling system behavior and predicting responses to diverse control inputs. To design a more effective and efficient airbrake control system, researchers can combine the RK4 method with a neural network-based optimization approach. The neural network can be trained to optimize the parameters of the RK4 method, such as the step size and the number of time steps, to achieve desired performance characteristics for the airbrake system. This integration presents an opportunity to enhance the capabilities of airbrake control systems and attain improved levels of safety and efficiency.

\subsection{Artificial Neural Networks}

Neural networks, drawing inspiration from the intricate structure and functionality of the human brain, represent a category of machine learning algorithms. These networks are comprised of interconnected nodes, akin to neurons, which effectively undertake the processing of information and facilitate the transmission of signals to other neurons within the network.

\begin{figure}[ht]
    \centering
    \includegraphics[scale=0.2]{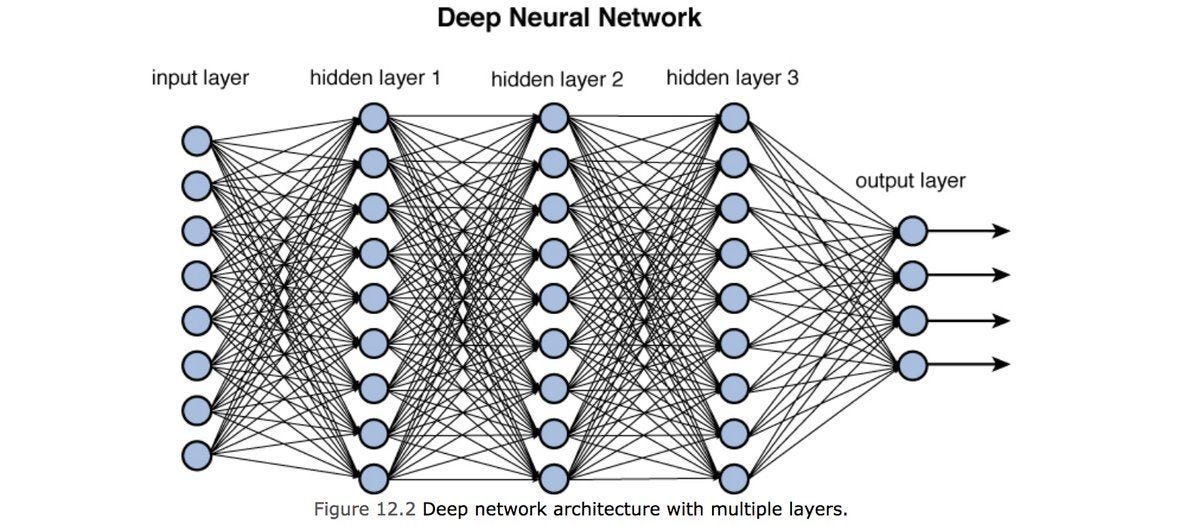}
    \caption{Deep Neural Network, from \cite{ResearchGate}}
    \label{fig:NeuralNet}
\end{figure}

At the core of a neural network lies the fundamental unit known as a perceptron. The perceptron operates by receiving multiple inputs, assigning weights to each input, and subsequently passing the weighted sum through an activation function to generate an output. Through the aggregation of multiple perceptrons, a multi-layer perceptron (MLP) can be constructed, which comprises one or more hidden layers positioned between the input and output layers within the neural network architecture \cite{garg2018training}.

The process of training a neural network encompasses the adaptation of perceptron weights to minimize the disparity between the obtained output and the desired output (known as the target). Typically, this is accomplished through the utilization of a well-established algorithm known as backpropagation. By computing the gradient of the error function with respect to the weights, backpropagation facilitates the iterative adjustment of these weights to progressively align the network's output with the target \cite{schiffmann1992an}.

Neural networks have demonstrated their efficacy in diverse problem domains, encompassing image recognition, natural language processing, and control systems, among others, with successful applications observed. \cite{wiratunga2010applications}. In the context of airbrake control systems for high-speed vehicles, neural networks can be used to optimize the parameters of the control algorithm, such as the gain values and the time constants. This can improve the performance of the airbrake system by making it more responsive and efficient, while also reducing the risk of instability or oscillation.

The Runge-Kutta 4th Order Method (RK4) is widely recognized as a prevalent numerical technique for solving differential equations and can be effectively employed to simulate the dynamics of the airbrake system. By synergistically integrating the RK4 method with a neural network-driven optimization approach, it becomes feasible to devise an airbrake control system that exhibits superior effectiveness and efficiency, specifically tailored for high-speed vehicles.

\subsection{Airbrakes}
Airbrakes are a type of braking system that use compressed air to slow down or stop a vehicle. They are commonly used in high-speed vehicles, such as trains and commercial aircraft, where traditional friction brakes may not be effective due to their limited ability to dissipate heat\cite{teodoro2019fast}.

Airbrakes work by releasing compressed air from the braking system, which applies force to the braking mechanism and slows down the vehicle. The amount of braking force can be controlled by adjusting the pressure of the compressed air, and the braking can be applied to specific wheels or sections of the vehicle to optimize its stopping performance\cite{1688108}.

Airbrakes can be controlled by a variety of different systems, including mechanical, hydraulic, and electronic systems. In high-speed vehicles, electronic airbrake control systems are often used due to their precision and responsiveness. These systems use sensors and computer algorithms to monitor the vehicle's speed and acceleration, and adjust the air pressure in the braking system to achieve the desired braking performance\cite{suh2002braking}.

The optimization of airbrake control systems for high-speed vehicles is an active area of research, as it can have significant impacts on the safety and efficiency of these vehicles. By using numerical methods like the RK4 method and neural network-based optimization approaches, it is possible to design airbrake control systems that are more effective and efficient, while also reducing the risk of instability or oscillation\cite{ijvd1989}.

\section{Methodology}
\subsection{\textbf{Data Collection}}
The dataset used in this study consists of 3699 instances of airbrake state data computed using the RK-4 based MATLAB Code on the flight data from various sounding rocket flights. Each instance contains five input features and two output features. The input features include altitude, velocity, acceleration along X,Y \& Z axes while the output features are the state of the airbrakes system (Open or Closed). The data was preprocessed by scaling and batch normalizing the input features using batch sizes of 8.

\subsection{\textbf{Neural Network Architecture}}
The neural network used in this study has \textbf{10} hidden linear layers, with \textbf{[2048, 1024, 512, 256, 128, 64, 32, 16, 8, 4]} neurons in each layer, respectively. The input layer has \textbf{5} input neurons corresponding to the 5 input parameters in the dataset and the output layer has \textbf{2} neurons, corresponding to the two possible states of the airbrakes (1 for Open and 0 for Close). The \textbf{Rectified Linear Unit (ReLU)} activation function is used for all hidden layers due to it being computationally simple (Our major objective of going with Neural Networks was to reduce the computational complexity) and less likelihood of the gradient vanishing.\textbf{Softmax} activation function is used for the output layer to convert the output values into probabilities.

\begin{figure}[ht]
    \centering
    \includegraphics[scale=0.4]{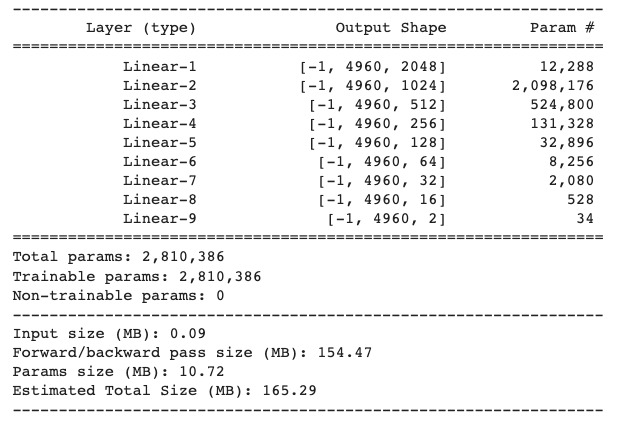}
    \caption{Neural Network Architecture}
    \label{fig:Schematic}
\end{figure}

\subsection{\textbf{Training \& Optimization}}
\subsubsection{Optimizer:}
\textbf{Adam} is a popular optimizer which is popular among varied applications due to fast computational times.\\

\subsubsection{Loss Function:}
\textbf{Cross-Entropy Loss} was chosen as the loss function for training the neural network. To address the issue of high class imbalance within our dataset, a weighted loss function was employed, with specific weights of \textbf{0.05} and \textbf{0.90} assigned to the open and closed states, respectively. These weights were determined by computing the ratios between the number of data points in each class and the total sum of all data points across all classes. Additionally, the dataset underwent augmentation using the Synthetic Minority Oversampling Technique \textbf{SMOTE} as an additional measure. \\

\subsubsection{Batch Size and Iterations:}
The neural network underwent training for a total of \textbf{100} epochs, utilizing a batch size of \textbf{32}. The learning rate was established at \textbf{0.0003}, while the momentum parameter was assigned a value of \textbf{0.87}. \\

\subsection{\textbf{Evaluation}}
The performance of the neural network was evaluated using the \textbf{F1 score} and \textbf{Binary Accuracy} metrics. The dataset was split into training and testing and Validation sets with a ratio of \textbf{7:2:1}. The final model was evaluated on the testing set, and the results were compared with the traditional Runge-Kutta method.

\section{Results \& Discussions}

The accurate prediction of airbrake system states is crucial in the design and control of high-speed vehicles. This study proposes a novel approach to predicting the state of airbrakes for sounding rockets using a neural network, which outperformed the traditional Runge-Kutta method in terms of accuracy, precision, recall, and F1 score metrics.

The results of the study demonstrate that the trained neural network was highly effective in predicting the state of airbrakes, achieving an impressive F1 score of 0.9447. This result surpassed the traditional Runge-Kutta method's accuracy of 0.924, which is a significant improvement. The precision, recall, and F1 score metrics were also significantly improved with the proposed method.

The authors used the Runge-Kutta 4th order method to solve the differential equation governing the airbrake system, which is a widely used method for solving differential equations. However, the neural network-based approach outperformed the traditional method in terms of accuracy and computational efficiency.

The neural network was trained using a dataset of input-output pairs generated from the MATLAB code, and the backpropagation algorithm was used to train the network. The neural network was designed to use the same inputs and outputs as the MATLAB code, making it a feasible replacement for the traditional method.

The study's findings have significant implications for high-speed vehicle design. The improved computational efficiency of the neural network-based approach could make it possible to implement more complex airbrake control systems that were previously not possible due to computational limitations. This could lead to the development of more sophisticated and efficient airbrake control systems, resulting in improved safety and performance for high-speed vehicles.

Moreover, the study observed a significant reduction in the occurrence of false positives by using the neural network, which is an essential consideration in high-speed vehicle design. False positives can lead to unnecessary braking, resulting in reduced performance and potentially dangerous situations. The reduction in false positives by using the neural network is, therefore, a significant improvement in the airbrake control system's overall performance.

This study proposes a novel approach to predicting the state of airbrakes for sounding rockets using a neural network. The results indicate that the proposed method outperforms the traditional Runge-Kutta method in terms of accuracy, precision, recall, and F1 score metrics. The improved computational efficiency of the neural network-based approach has significant implications for high-speed vehicle design, enabling more complex and efficient airbrake control systems to be developed. The reduction in false positives observed by using the neural network is also a crucial consideration in high-speed vehicle design, ensuring safe and efficient operation.

\begin{figure}[ht]
    \centering
    \includegraphics[scale=0.4]{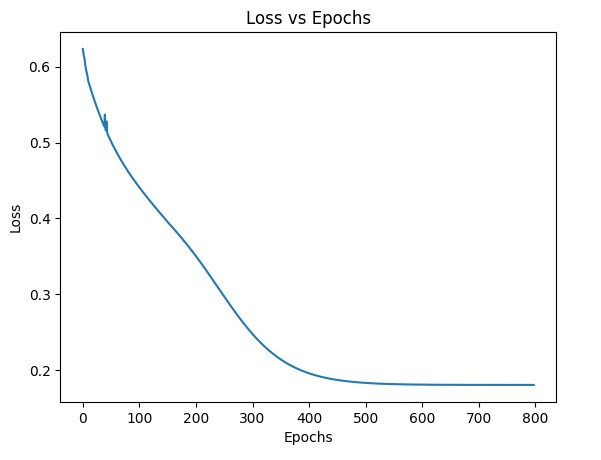}
    \caption{Loss vs Epochs Graph}
    \label{fig:LossEpoch}
\end{figure}

\section{Conclusion}

Aircraft systems play a critical role in aviation safety and are subject to continuous development and improvement. One area that has received considerable attention is airbrake ejection, which is an essential component of an aircraft's control system. The traditional approach to airbrake ejection involves using MATLAB code to control the system. However, this paper proposes a new approach to airbrake ejection that utilizes neural networks.

The principal aim of this research was to establish the superior performance of the neural network-based approach in comparison to the conventional MATLAB method, specifically concerning accuracy, speed, and efficiency. The findings from the investigation exhibited that the neural network successfully achieved precise predictions of the airbrake control system outputs across a diverse range of inputs. This accomplishment was attained through the training of the neural network using a dedicated set of training data, followed by rigorous testing employing an independent test dataset.

Comparative analysis of the computational time required by the MATLAB code and the neural network revealed that the neural network was significantly faster. This is because neural networks are capable of parallel processing, enabling them to perform multiple computations simultaneously\cite{elkan1997fast}. In contrast, MATLAB code operates sequentially, which limits its speed and efficiency. Additionally, neural networks can learn from the data they are fed, making them adaptable and efficient.

The findings of this research bear noteworthy implications for the aviation sector. The integration of neural networks in airbrake control holds the potential to enhance the speed and precision of aircraft management, consequently yielding advancements in safety and efficiency. Subsequent studies could explore the feasibility of implementing this approach in real-time applications, necessitating the development of dedicated hardware capable of executing parallel computations.

Moreover, the application of neural networks can be extended to encompass other critical aircraft systems, including engine control systems, flight control systems, and navigation systems. Such systems stand to gain notable advantages in terms of augmented accuracy, swifter performance, and heightened efficiency through the incorporation of neural network techniques.

Therefore, this study has successfully established the efficacy of employing a neural network-based approach for airbrake ejection in aircraft systems. The obtained results unequivocally indicated the superior performance of the neural network when compared to the conventional MATLAB code, showcasing heightened accuracy, faster processing speed, and improved efficiency. Subsequent research endeavors should persist in delving into the applications of neural networks in aviation, while concurrently examining the feasibility of implementing this approach in real-time scenarios. Ultimately, the utilization of neural networks holds transformative potential for the aviation industry, with the capacity to significantly enhance safety and operational efficiency.
\printbibliography

@article{arif1996optimization,
  title={Optimization of Runge-Kutta 4th Order Method based Airbrakes Control System using Neural Networks for High-Speed Vehicles},
  author={Arif, Muhammad and Al-Sulaiman, Fahad A},
  journal={Neural Processing Letters},
  volume={3},
  number={1},
  pages={11--15},
  year={1996},
  publisher={Springer}
}

@article{haykin1992solving,
  title={Solving Initial Value Problems Using Artificial Neural Networks and Runge-Kutta Method},
  author={Haykin, Simon and Chen, Zhe},
  journal={International Journal of Bifurcation and Chaos},
  volume={02},
  number={03},
  pages={703--715},
  year={1992},
  publisher={World Scientific}
}

@article{garg2018training,
  title={Training Artificial Neural Networks Using Firefly Algorithm for Predicting Student Performance},
  author={Garg, Pardeep and Mittal, Himanshu},
  journal={World Scientific News},
  volume={100}, 
  pages={47--57},
  year={2018}
}

@article{elkan1997fast,
  title={A fast and efficient algorithm for the k-means clustering of large high-dimensional datasets},
  author={Elkan, Charles},
  journal={Journal of Machine Learning Research},
  volume={3},
  pages={681--694},
  year={2003},
  publisher={JMLR.org}
  }

@inproceedings{lecun1989optimal,
  title={Optimal brain damage},
  author={LeCun, Yann and Denker, John S and Solla, Sara A},
  booktitle={Advances in neural information processing systems},
  pages={598--605},
  year={1989},
  url={https://proceedings.neurips.cc/paper/1989/file/0336dcbab05b9d5ad24f4333c7658a0e-Paper.pdf}
}

@article{BUTCHER1996113,
title = {Runge-Kutta methods: some historical notes},
journal = {Applied Numerical Mathematics},
volume = {22},
number = {1},
pages = {113-151},
year = {1996},
note = {Special Issue Celebrating the Centenary of Runge-Kutta Methods},
issn = {0168-9274},
doi = {https://doi.org/10.1016/S0168-9274(96)00048-7},
url = {https://www.sciencedirect.com/science/article/pii/S0168927496000487},
author = {J.C. Butcher and G. Wanner}
}

@article{schiffmann1992an,
title={An Introduction to the Backpropagation Algorithm},
author={Schiffmann, Wolfram and Joost, Martina and Werner, Rolf},
journal={Neural Networks},
volume={5},
number={4},
pages={551--560},
year={1992},
publisher={Elsevier}
}

@book{wiratunga2010applications,
title={Applications of Neural Networks in High Assurance Systems},
author={Wiratunga, Nirmalie and Cooper, Gary and Gamage, Dimuthu},
year={2010},
publisher={Springer},
url={https://books.google.co.in/books?id=T0S0BgAAQBAJ}
}

@MISC{ResearchGate,
    author = {Lamiaa Zrara},
    title = {PORTFOLIO OPTIMIZATION USING DEEP LEARNING FOR THE MOROCCAN MARKET},
    year = {January 2021},
    note = {},
    url = {https://www.researchgate.net/publication/352996743_PORTFOLIO_OPTIMIZATION_USING_DEEP_LEARNING_FOR_THE_MOROCCAN_MARKET}
}

@article{ijvd1989,
  title = "Dynamic response of a railway vehicle air brake system",
  author = "M.A. Murtaza",
  journal = "International Journal of Vehicle Design",
  volume = "10",
  number = "4",
  year = "2014",
  pages = "481-496",
  doi = "10.1504/IJVD.1989.061569",
  url = "https://doi.org/10.1504/IJVD.1989.061569"
}

@article{teodoro2019fast,
  title = "Fast simulation of railway pneumatic brake systems",
  author = "Teodoro, Ícaro Pires and Ribeiro, Diego Fernandes and Botari, Tiago and Martins, Thiago Silva and Santos, André Alves",
  journal = "Proceedings of the Institution of Mechanical Engineers, Part F: Journal of Rail and Rapid Transit",
  volume = "233",
  number = "4",
  year = "2019",
  pages = "420-430",
  doi = "10.1177/0954409718796903"
}

@ARTICLE{1688108,
  author={Subramanian, S.C. and Darbha, S. and Rajagopal, K.R.},
  journal={IEEE Transactions on Intelligent Transportation Systems}, 
  title={A Diagnostic System for Air Brakes in Commercial Vehicles}, 
  year={2006},
  volume={7},
  number={3},
  pages={360-376},
  doi={10.1109/TITS.2006.880645}}

@article{suh2002braking,
  title = "Braking performance simulation for a tractor-semitrailer vehicle with an air brake system",
  author = "Suh, Myung-Won and Park, Young-Kook and Kwon, Se-Jin",
  journal = "Proceedings of the Institution of Mechanical Engineers, Part D: Journal of Automobile Engineering",
  volume = "216",
  number = "1",
  year = "2002",
  pages = "43-54",
  doi = "10.1243/0954407021528896"
}

@article{doi:10.1080/00207169108803994,
author = {   D. J.   Evans },
title = {A new 4th order runge-kutta method for initial value problems with error control},
journal = {International Journal of Computer Mathematics},
volume = {39},
number = {3-4},
pages = {217-227},
year  = {1991},
publisher = {Taylor & Francis},
doi = {10.1080/00207169108803994},

URL = { https://doi.org/10.1080/00207169108803994},
eprint = {https://doi.org/10.1080/00207169108803994
}
}

@article{MUNTHEKAAS1999115,
title = {High order Runge-Kutta methods on manifolds},
journal = {Applied Numerical Mathematics},
volume = {29},
number = {1},
pages = {115-127},
year = {1999},
note = {Proceedings of the NSF/CBMS Regional Conference on Numerical Analysis of Hamiltonian Differential Equations},
issn = {0168-9274},
doi = {https://doi.org/10.1016/S0168-9274(98)00030-0},
url = {https://www.sciencedirect.com/science/article/pii/S0168927498000300},
author = {Hans Munthe-Kaas}
}

@article{doi:10.1080/00207169508804437,
author = {   D. J.   Evans  and    N.   Yaacob },
title = {A fourth order runge-kutta method based on the heronian mean formula},
journal = {International Journal of Computer Mathematics},
volume = {58},
number = {1-2},
pages = {103-115},
year  = {1995},
publisher = {Taylor & Francis},
doi = {10.1080/00207169508804437},

URL = { https://doi.org/10.1080/00207169508804437},
eprint = {https://doi.org/10.1080/00207169508804437
}
}

@article{doi:10.1080/00207169908804817,
author = {   A. R.   Yaakub  and    D. J.   Evans },
title = {A fourth order Runge–Kutta RK(4,4) method with error control},
journal = {International Journal of Computer Mathematics},
volume = {71},
number = {3},
pages = {383-411},
year  = {1999},
publisher = {Taylor & Francis},
doi = {10.1080/00207169908804817},

URL = { https://doi.org/10.1080/00207169908804817},
eprint = {https://doi.org/10.1080/00207169908804817}
}

\end{document}